\documentclass[11pt]{article}

\usepackage[latin1]{inputenc}

\usepackage{amsfonts,amsmath,amssymb,amsthm,graphicx,epsfig,float}
\usepackage[T1]{fontenc}

\usepackage[english,francais]{babel}

{
  
  \newtheorem*{theorem*}{Theorem}
  \newtheorem{lemma}{Lemma}
  \newtheorem{corollary}{Corollary}

\newtheorem*{corollaire*}{Corollaire}
\newtheorem*{proposition*}{Proposition}
\theoremstyle{remark}
  \newtheorem*{remarque*}{Remarque}
}

\newcounter{ex}

\newenvironment{rem*}{
  \noindent\textbf{Remarque. }}{}

\newcommand{\Nn}{\mathbb{N}}

\newcommand{\Scal}{\mathcal{S}}
\renewcommand {\epsilon}{\varepsilon}

\renewcommand {\leq}{\leqslant}
\renewcommand {\geq}{\geqslant}

\title{{\bf Construction of measures with dilation}}
\author{Henry de Thélin}
\date{}

\begin{document}
\maketitle

%$\Cc, \Rr$%

%% Redefinition Titre
\def\figurename{{Fig.}}%
\def\proofname{Preuve}% for AMS-\LaTeX
\def\contentsname{Sommaire}%
%% Fin

\selectlanguage{english}

\begin{abstract}

We give a construction of measures with partial sum of Lyapunov
exponents bounded by below.

\end{abstract}

AMS: 28Dxx, 58F11.

\section*{{\bf Introduction}}
\par

Let $M$ be a compact $C^1$-Riemannian manifold of dimension $d$ and let $f: M
\mapsto M$ be a $C^1$-map.

For $1 \leq k \leq d$, we denote by $\Scal_k$ the set of $C^1$-maps $\sigma :D^k= [0,1]^k \mapsto M$. We define the $k$-volume of $\sigma \in \Scal_k$ with the formula:

$$ V(\sigma)=\int_{D^k} | \Lambda^k T_x \sigma | d \lambda(x),$$
where $d \lambda$ is the Lebesgue measure on $D^k$ and  $|\Lambda^k T_x \sigma|$ is the
norm of the linear map $\Lambda^k T_x \sigma : \Lambda^k T_x D^k
\mapsto \Lambda^k T_{\sigma(x)} M$ induced by the Riemannian metric on
$M$.

Some links between the volume growth of iterates of submanifolds of
$M$ and the entropy of $f$ have been studied by Y. Yomdin (see
\cite{Y} and \cite{Gr}), S. E. Newhouse (see
\cite{Ne}), O.S. Kozlovski (see \cite{Ko}) and J. Buzzi (see \cite{Bu}).

In this article, we prove that the volume growth of iterates of submanifolds
of $M$ permits to create invariant measures with partial sum of Lyapunov
exponents bounded by below. More precisely, for $1 \leq k \leq d$ we define the $k$-dilation:

$$d_k:= \limsup_{n \rightarrow \infty} \frac{1}{n} \log \sup_{\sigma \in \Scal_k}
\frac{V(f^n \circ \sigma)}{V(\sigma)}.$$

We will prove the following theorem:

\begin{theorem*}

For all integer $k$ between $1$ and $d=\mbox{dim}(M)$ there exists an ergodic measure $\nu(k)$ for
which:
$$ \sum_{i=1}^{k} \chi_i \geq d_k.$$  

Here $\chi_1 \geq \chi_2 \geq \dots \geq \chi_d$ are the Lyapunov
exponents of $\nu(k)$.
\end{theorem*}

Notice that when $k=d$ and $f$ is a ramified covering in some sense, the theorem can be deduced from a result due
to T.-C. Dinh and N. Sibony (see \cite{DS} paragraph 2.3).

\section*{{\bf Proof of the theorem}}

Let $k$ be a positive integer between $1$ and $d$. We have to prove that there
exists an ergodic measure $\nu(k)$ for which 
$$\sum_{i=1}^{k} \chi_i = \lim_{m \rightarrow \infty}
\frac{1}{m} \int \log | \Lambda^k T_y f^m | d \nu(k)(y) \geq d_k.$$

For the definition of Lyapunov exponents
and for the statement of the previous equality, see \cite{KH} and
\cite{Ar} chapter $3$. 

$ $

There will be three steps in the proof of the theorem.

In the first one, we will change the dilation $d_k$ into a dilation of
$ |\Lambda^k T_{x} f^{n}|$. More precisely, we will find points
$x_{n_l}$ with $\frac{1}{n_l}
\log |\Lambda^k T_{x(n_l)} f^{n_l}| \geq d_k - \epsilon$.

In the second part, we will see that the dilation of $|\Lambda^k
T_{x(n_l)} f^{n_l}|$ can be spread out in time. We will give the construction
of a measure $\nu_l$ such that $d_k- 2 \epsilon \leq \frac{1}{m} \int \log | \Lambda^k T_y f^m | d
  \nu_l(y)$.

The third step of the proof will be to take the limit in the previous
inequality.

\subsection*{{\bf 1) First step}}

Let $n_l$ be a subsequence such that:

$$ \frac{1}{n_l} \log \sup_{\sigma \in \Scal_k}
\frac{V(f^{n_l} \circ \sigma)}{V(\sigma)}  \rightarrow  d_k.$$

We can find now a sequence $\sigma_{n_l} \in \Scal_k$ which verifies:

$$ \frac{1}{n_l} \log \frac{V(f^{n_l} \circ
  \sigma_{n_l})}{V(\sigma_{n_l})}  \rightarrow  d_k.$$

In the next lemma, we prove that we have dilation for $| \Lambda^k T_x
f^n|$ for some $x$:

\begin{lemma}

For all $l \geq 0$ there exists $x(n_l) \in M$ with:

$$ \log |\Lambda^k T_{x(n_l)} f^{n_l}| \geq \log \left( \frac{ V(f^{n_l}
\circ  \sigma_{n_l})}{2V(\sigma_{n_l})} \right).$$

\end{lemma}

\begin{proof}

Otherwise we would have an integer $l$ such that for all $x \in M$:

$$  |\Lambda^k T_{x} f^{n_l}| \leq  \frac{ V(f^{n_l}
\circ  \sigma_{n_l})}{2V(\sigma_{n_l})} .$$

So (see \cite{Ar} chapter 3.2.3 for properties on exterior powers),

$$ V(f^{n_l} \circ  \sigma_{n_l})= \int_{D^k} |\Lambda^k T_x(f^{n_l} \circ
\sigma_{n_l})| d \lambda(x) =   \int_{D^k} |\Lambda^k
T_{\sigma_{n_l}(x)}f^{n_l} \circ \Lambda^k T_x \sigma_{n_l}| d \lambda(x)$$
is bounded by above by

$$ \int_{D^k} |\Lambda^k
T_{\sigma_{n_l}(x)}f^{n_l} | | \Lambda^k T_x \sigma_{n_l})| d \lambda(x)
\leq \frac{ V(f^{n_l} \circ  \sigma_{n_l})}{2}$$

and we obtain a contradiction.

\end{proof}

\begin{corollary}

There exists a sequence $\epsilon(l)$ which converges to $0$ such that:

$$\frac{1}{n_l}  \log |\Lambda^k T_{x(n_l)} f^{n_l}| \geq d_k - \epsilon(l),$$
for some points $x(n_l)$ in $M$.

\end{corollary}

\subsection*{{\bf 2) Second step}}

In this section, we will spread out in time the previous dilation.

Let $m$ be a positive integer. We will now cut $n_l$ with $m$ different ways.

By using the Euclidian division, we can find $q_l^i$ and $r_l^i$ (for $i=0, \dots ,m-1$) such that:

$$n_l=i+ m \times q_l^i + r_l^i$$
with $0 \leq r_l^i < m$.

If $i \in \{0, \dots , m-1 \}$, we have:

$$| \Lambda^k T_{x(n_l)} f^{n_l} | \leq  |\Lambda^k T_{f^{i+ mq_l^i}(x(n_l))}f^{r_l^i}| \times
\prod_{j=0}^{q_l^i-1} | \Lambda^k T_{f^{i+jm}(x(n_l))} f^m | \times |
\Lambda^k T_{x(n_l)} f^i |,$$

so, by using the previous corollary,

$$n_l(d_k- \epsilon(l)) \leq  \log |\Lambda^k T_{f^{i+
    mq_l^i}(x(n_l))}f^{r_l^i}| + \sum_{j=0}^{q_l^i-1} \log | \Lambda^k
    T_{f^{i+jm}(x(n_l))} f^m | + \log |
\Lambda^k T_{x(n_l)} f^i |.$$

If we take the sum on the $m$ different ways to write $n_l$, we
obtain:

$$ m n_l (d_k - \epsilon(l)) \leq \sum_{i=0}^{m-1} \log |\Lambda^k T_{f^{i+
    mq_l^i}(x(n_l))}f^{r_l^i}|+\sum_{i=0}^{m-1} \sum_{j=0}^{q_l^i-1} \log | \Lambda^k
    T_{f^{i+jm}(x(n_l))} f^m | + \sum_{i=0}^{m-1} \log |
\Lambda^k T_{x(n_l)} f^i |.$$

We have to transform this estimate on a relation on a measure. To
realize that, we remark that:

$$\log | \Lambda^k T_{f^p(x(n_l))} f^m | = \int \log | \Lambda^k T_y
f^m | d \delta_{f^p(x(n_l))}(y),$$
where $\delta_{f^p(x(n_l))}$ is the dirac measure at the point
$f^p(x(n_l))$.

So the previous inequality becomes:

$$d_k- \epsilon(l) \leq a_l + \frac{1}{m} \int \log | \Lambda^k T_y
f^m | d \left(
\frac{1}{n_l} \sum_{i=0}^{m-1} \sum_{j=0}^{q_l^i-1}
\delta_{f^{i+mj}(x(n_l))} \right) (y) + b_l$$

with
$$a_l= \frac{1}{mn_l} \sum_{i=0}^{m-1} \log |\Lambda^k T_{f^{i+
    mq_l^i}(x(n_l))}f^{r_l^i}| $$
and 

$$b_l =\frac{1}{m n_l}  \sum_{i=0}^{m-1} \log | \Lambda^k T_{x(n_l)} f^i
|.$$

Now, because $f$ is a $C^1$-map we have:
$$a_l \leq \frac{1}{m n_l} \sum_{i=0}^{m-1} \log L^{mk} \leq \frac{k
  m^2}{mn_l} \log L$$
where $L = \max(\max_x |T_x f|,1)$ and:
$$b_l \leq \frac{1}{m n_l} \sum_{i=0}^{m-1} \log L^{mk} \leq \frac{k
  m^2}{m n_l} \log L.$$
So the sequences $a_l$ and $b_l$ are bounded by above by a sequence
  which converges to $0$ when $l$ goes to infinity.

In conclusion, we have:

\begin{equation}{\label{eq1}}
d_k-  \epsilon'(l) \leq \frac{1}{m} \int \log | \Lambda^k T_y f^m | d
  \nu_l(y)
\end{equation}
with 
$$ \nu_l=\frac{1}{n_l} \sum_{i=0}^{m-1} \sum_{j=0}^{q_l^i-1}
\delta_{f^{i+mj}(x(n_l))},$$
and $\epsilon'(l)$ a sequence which converges to $0$.

\subsection*{{\bf 3) Third step}}

The aim of this section is to take a limit for $\nu_l$ in the
equation (\ref{eq1}).

First, observe that $\nu_l= \frac{1}{n_l} \sum_{p=0}^{n_l-m}
\delta_{f^p(x(n_l))}$ and that the sequence $ \frac{1}{n_l}
\sum_{p=0}^{n_l-1} \delta_{f^p(x(n_l))} - \nu_l$ converges to $0$. In
particular, there exists a subsequence of $\nu_l$ which converges
to a measure $\nu$ which is a probability invariant under $f$ and independant of $m$. We continue to call $\nu_l$ the subsequence which converges
to $\nu$. To complete the proof of the theorem,
we have to take the limit in the equation (\ref{eq1}). However, we have
to be careful because the function $y \mapsto \log | \Lambda^k T_y f^m |$
is not continuous. But, we have the following lemma:

\begin{lemma}

$$\limsup_{l \rightarrow \infty} \frac{1}{m} \int \log | \Lambda^k T_y f^m | d \nu_l(y) \leq
\frac{1}{m} \int \log |\Lambda^k T_y f^m | d \nu(y).$$

\end{lemma}

\begin{proof}

For $r \in \Nn$, let $\Phi_r(y)= \max (\log | \Lambda^k T_y f^m |, -r)$. 

The functions $\Phi_r$ are continuous and the sequence $\Phi_r$ decreases to the map $y
\mapsto \log | \Lambda^k T_y f^m |$ when $r$ goes to infinity.

Then:

$$\frac{1}{m} \int \log | \Lambda^k T_y f^m | d \nu_l(y) \leq \frac{1}{m}
\int \Phi_r(y) d \nu_l(y),$$

and,

$$\limsup_{l \rightarrow \infty}  \frac{1}{m} \int \log | \Lambda^k T_y f^m | d \nu_l(y)
\leq \frac{1}{m} \int  \Phi_r(y) d \nu(y)$$

because $\Phi_r$ is continuous. Now, we obtain the lemma by using the monotone convergence
theorem.

\end{proof}

It remains to take the limit in the equation (\ref{eq1}). We obtain then
the

\begin{corollary}

For all $m$, we have:

$$ d_k \leq \frac{1}{m} \int \log |\Lambda^k T_y f^m | d \nu(y).$$

\end{corollary}

In particular,

$$ d_k \leq  \int \sum_{i=1}^{k} \chi_i(y) d \nu(y)$$

where the $\chi_1 \geq \chi_2 \geq \dots \geq \chi_d$ are the
Lyapunov exponents of $\nu$. Finally, by using the ergodic decomposition of $\nu$, we obtain the existence of an
ergodic measure $\nu(k)$ with:

$$ d_k \leq  \sum_{i=1}^{k} \chi_i.$$

\bigskip

Henry de Thélin\\
Université Paris-Sud (Paris 11)\\
Mathématique, Bât. 425\\
91405 Orsay\\
France

\end{document}